\documentclass[12pt]{article}
\usepackage{graphicx}
\usepackage{amsfonts}
\usepackage{amsmath}
\usepackage{amssymb}  
\usepackage{gastex}
\usepackage{color}

\def\R{{\mathbb{R}}}

\newcommand{\dst}{\displaystyle}

\def\be{\begin{equation}}
\def\ee{\end{equation}}
\def\ba{\begin{array}}
\def\ea{\end{array}}
\def\eqa{\begin{eqnarray}}
\def\eqe{\end{eqnarray}}

\newtheorem{corollary}{Corollary}
\newtheorem{proposition}{Proposition}

\newtheorem{lemma}{Lemma}

\newenvironment{definition}{\medskip\noindent{\it Definition. }}{ \medskip}
\newenvironment{proof}{\medskip\noindent{\it Proof.}}{\medskip}

\newenvironment{remark}{\medskip\noindent{\it Remark. }}{
\medskip}

\begin{document}

\title{Robust Stabilization of Nonlinear Systems by Quantized and Ternary Control
\thanks{C.~De Persis is
with the Dipartimento di Informatica e Sistemistica A.\ Ruberti, Sapienza Universit\`a
di Roma, Via Ariosto 25, 00185 Roma, Italy.
Fax +390677274033, Ph.~+390677274060, Email {\tt depersis@dis.uniroma1.it}.}}
\author{Claudio De Persis}

\date{September 26, 2008}

\maketitle

\begin{abstract}
Results on the problem of stabilizing a nonlinear continuous-time system by a
finite number of control or measurement values are presented. The
basic tool is a discontinuous version of the so-called semi-global
backstepping lemma. We derive robust practical stabilizability
results by quantized and ternary controllers and apply them to
some significant control problems.
\end{abstract}

\medskip


\section{Introduction}
The problem of controlling systems through a limited bandwidth channel
has recently raised a great interest in the community, as thoroughly
surveyed in \cite{nair.et.al.proc.ieee07}. A possible approach
to the problem for continuous-time systems consists of partitioning
(a subset of) the state space into a finite number of regions and
transmitting information whenever the state crosses one of the
boundaries. The resulting system is known as a quantized control
system, and the focus of this paper is on this class of systems.
Many authors have
contributed to the topic, and we refer the interested reader to
\cite{nair.et.al.proc.ieee07} for an exhaustive bibliography. Among
the papers which are important to our derivations we recall
\cite{liberzon.aut03}, \cite{elia.mitter.tac01}, \cite{hayakawa.et.al.acc06} and
\cite{ceragioli.depersis.scl07}.
\\
In \cite{liberzon.aut03},
adopting a time-varying quantization, and relying on input-to-state
stability of the system, the author shows asymptotic
convergence to the origin.
In \cite{hayakawa.et.al.acc06} and \cite{ceragioli.depersis.scl07},
the role of static  logarithmic
quantization (\cite{elia.mitter.tac01}) to prove practical semi-global
stabilizability of nonlinear stabilizable
systems has been
 investigated. The two papers mainly differ
in the type of solution adopted. In particular, the paper
\cite{ceragioli.depersis.scl07} establishes a few connections between quantized control and
discontinuous control systems, investigating Carath\'eodory and
Krasowskii solutions in the context of quantized control systems,
while the authors of \cite{hayakawa.et.al.acc06} propose a
hysteresis-based implementation of the quantized control.
The two papers also present results which rely on notions of robustness
different from input-to-state stability.
Finally in \cite{hayakawa.et.al.acc06},
an adaptive control scheme for nonlinear continuous time uncertain systems is
proposed.
\\
In this paper, we establish a few results on the problem of
stabilizing a nonlinear continuous-time  by
quantized control robustly with respect to uncertainties.
A discontinuous version of the {\em semi-global backstepping} lemma of
\cite{teel.praly.sicon95}, in which
the measured state is
logarithmically quantized, is applied to show that
{\em minimum-phase} nonlinear systems,
possibly with {\em uncertain parameters},
can be robustly semi-globally practically stabilized by a
quantized function of
{\em partial-state} measurements. The control techniques introduced
in the papers previously discussed can not be directly applied to
the problem considered here, and the resulting quantized control
we propose is new to the best of our knowledge. In the scenario in which the
feedback information travels through a finite-bandwidth channel, it
is important to calculate the bandwidth needed to
implement the quantized controller.
The solution we propose has the additional advantage of allowing us
to estimate an upper bound on the required bandwidth.
We also show that
semi-global practical stabilization is possible
even using a simple {\em switched ternary} controller.
The backstepping lemma of \cite{teel.praly.sicon95} has played a
fundamental role in the design of many robust nonlinear control
schemes (\cite{isidori.book.vol2}). We conclude the paper presenting a few examples where the quantized backstepping lemma
is used to solve some of these robust nonlinear control problems in the presence of
quantization.
\\
Preliminary facts are presented in Section \ref{sec.prel}.
In Section \ref{sec.quant}, the semi-global backstepping tool in the
presence of quantization is proven. An upper  bound on the
bandwidth associated with the quantized control scheme we propose is
estimated in Section \ref{sec.band}.
The ternary controller is
introduced in Section \ref{sec.ternary}. Some examples are illustrated in Section
\ref{sec.appl}.

\vspace{-0.25cm}

\section{Preliminaries}\label{sec.prel}
The system we focus our attention on is
of the form
\vspace{-0.25cm}
\be\label{basic}\ba{rcl}
\dot x &=& F(x,\mu)+G(x,\mu)\zeta\\
\dot \zeta &=& q(x,\zeta,\mu)+b(x,\zeta,\mu)u
\ea\ee
with $x\in \mathbb{R}^{n-1}$, $\zeta\in \mathbb{R}$, $\mu$
an unknown parameter ranging over the compact set $\mathcal{
P}$, $u\in
\mathbb{R}$, $b(x,\zeta,\mu)\ge b_0>0$ for all $(x,\zeta,\mu)$.
The role of this kind of system to solve many important control
problems will be emphasized later on (cf.\ Section \ref{sec.appl}).
We suppose that the upper subsystem
satisfies the following property (\cite{teel.praly.sicon95}, see also
\cite{isidori.book.vol2}), which claims that the upper subsystem with $\zeta=0$
is asymptotically stable with a given region of attraction:

\begin{definition}\label{ulp}
The system
$\dot x=F(x,\mu)$,  $x\in \mathbb{R}^{n-1}$,
satisfies a Uniform Lyapunov Property if there exists an open set
$\mathcal{A}\subset \mathbb{R}^{n-1}$, a real number $c\ge 1$, a
continuously differentiable definite positive function $V\,:\,\mathcal{A}\to \R_+$ such
that $\Gamma_{c+1}:=\{x\,:\, V(x)\le c+1\}\subset \mathcal{A}$ and
$
\frac{\partial V}{\partial x} F(x,\mu) <0
$, for all $x\in
\Gamma_{c+1}$, $x\ne 0$,  for all  $\mu\in\mathcal{P}$.  $\triangleleft$
\end{definition}

\begin{figure}
\begin{center}
\scalebox{0.25}{\includegraphics{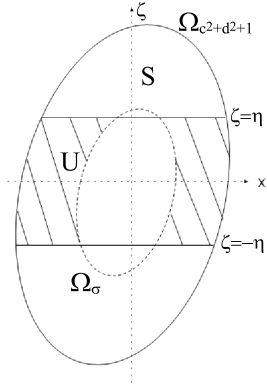}}
\caption{\label{fig.sets}
The figure represents the sets of
interest in the paper. The outer contour is
the boundary of $\Omega_{c^2+d^2+1}$, while the
inner contour is the boundary of $\Omega_\sigma$. $S$ is the region
between the 2. The 2 horizontal segments
are the set of points such that $\zeta=\pm \eta$.
The open set $U$ is emphasized by oblique solid lines.
The regions at the top, center and
bottom, delimited by the boundary of $\Omega_{c^2+d^2+1}$ and the 2
horizontal solid lines,
are respectively $\Omega_-$, $\Omega_0$, $\Omega_+$.
}
\end{center}
\end{figure}

Introduce the Lyapunov function (\cite{teel.praly.sicon95})
\[
W(x,\zeta)=\dst\frac{c V(x)}{c+1-V(x)} +
\dst\frac{d \zeta^2}{d+1-\zeta^2}
\]
defined on the
set $\{x \,:\, V(x)< c+1\}\times \{\zeta \,:\, \zeta^2< d+1\}$,
for some $d\ge 1$, and definite positive and proper therein.
For an arbitrary $\sigma>0$,
consider the set
$
S=
\{
(x,\zeta)\,:\, \sigma\le W(x,\zeta) \le c^2+d^2+1
\}
$.
The set is well-defined, because 
if $W(x,\zeta) \le c^2+d^2+1$, then $V(x)<c+1$ and $\zeta^2<d+1$.
In \cite{teel.praly.sicon95} (see also \cite{bacciotti.tac89}) it is proven
that a linear high-gain {\em partial-state} feedback $u=\bar
k \zeta$ exists which makes $\dot W(x,\zeta)$ negative on $S$
(thus allowing the authors to conclude that any trajectory
starting in $S$ is attracted by $\Omega_\sigma:=\{(x,\zeta)\,:\,
W(x,\zeta) \le \sigma\}$). In this paper, we are interested to
carry out an analogous investigation in the cases in which the feedback information
$\zeta$ is available in a ``limited" form, namely it undergoes
quantization.
\\
Following \cite{teel.praly.sicon95}, consider
the derivative $\dot W(x,\zeta)=(\partial W/\partial
x) \dot x+(\partial W/\partial
\zeta) \dot \zeta$.
It is possible to obtain the following inequality to hold
for all $(x,\zeta)\in
\Omega_{c^2+d^2+1}$:
\be\label{w.dot}
\dot W(x,\zeta)\le
\dst\frac{c}{c+1}\dst\frac{\partial V}{\partial
x}F(x,\mu)
+w(x,\zeta,\mu)\zeta+ 2\dst\frac{d (d+1)}{(d+1-\zeta^2)^2}\zeta b(x,\zeta,\mu)u
\;,
\ee
where
$w(x,\zeta,\mu)=\dst\frac{c(c+1)}{(c+1-V(x))^2}
\dst\frac{\partial V}{\partial
x}G(x,\mu)+2\dst\frac{d(d+1)}{(d+1-\zeta^2)^2}q(x,\zeta,\mu)$.\\
Because of the ULP property, if the state belongs to
$S_0=
\{
(x,\zeta)\in S\,:\, \zeta=0
\}$,
then $\dot W(x,\zeta)<0$.
By continuity, there exists a
neighborhood $U$ of $S_0$ where
the sum of the first two terms on the right-hand side of the inequality above
remains
strictly negative. Without loss of generality, we can suppose that a
constant $\eta>0$ exists such that $U=\{(x,\zeta)\in S\,:\,
|\zeta|<\eta\}$ (see Figure \ref{fig.sets}). Then,
to show that  $\dot W(x,\zeta)$ is negative on $S$,
it is enough to investigate the sign of
$\dot W(x,\zeta)$ on  $\tilde S:=S\setminus U$ only.

\section{Stabilization by quantized control}\label{sec.quant}

In what follows,
we consider the case in which the control $\bar k \zeta$, or the measurement $\zeta$,
is quantized by a logarithmic quantizer.
Let $u_0\in \R_+$, $j\in \mathbb{N}$ and $0<\delta<1$ be constants to design. Also let
$u_i=\rho^i u_0$,
with $\rho=\frac{1-\delta}{1+\delta}$ (\cite{elia.mitter.tac01}) . The following map is the
quantizer
\be\label{ppsi}
\Psi(r)=
\left\{\ba{crclrcll}
u_i &\dst\frac{1}{1+\delta}u_i& < & r &\le &
\dst\frac{1}{1-\delta}u_i\;, & 0\le i\le j \\[2mm]
0 & 0&\le & r &\le & \dst\frac{1}{1+\delta}u_j & \\[2mm]
-\Psi(-r) & & & r & < & \hspace{-1cm} 0\;, &
\ea\right.
\ee
and $u=-\Psi(\bar k \zeta)$ is the quantized input. We do not
define the quantizer for $\bar k \zeta>(1-\delta)^{-1}u_0$, since
$u_0$ will be designed in such a way that this bound is never
exceeded.
Observe that it
is {\em equivalent} to consider either the quantized control law $u=-\Psi(\bar k
\zeta)$ or the control law $u=-\bar k \bar \Psi(\zeta)$,  function
of the
{\em quantized partial-state} $\bar \Psi(\zeta)$, provided that $\bar \Psi$
is appropriately defined.
As a matter of fact, define $\bar
\Psi$ as $\Psi$ in (\ref{ppsi}), but with a new set of
quantization levels $\bar u_i$ (instead of $u_i$) defined as $\bar u_i=\rho^i
\bar
u_0$, with $\bar u_0=\bar k^{-1} u_0$. Then, it is easy to show
that $\bar k \bar \Psi(\zeta)=\Psi(\bar k \zeta)$, and all the results drawn
with $u=-\Psi(\bar k \zeta)$ also hold for $u=-\bar k
\bar\Psi(\zeta)$. In what follows, we only refer to the quantized input
$u=-\Psi(\bar k \zeta)$.
\\
Observe that
the quantizer has $2j+3$ quantization levels,
with $u_0$, $j$, $\bar k$ to determine. Of course, the size
of the {\em deadzone} of the quantizer, i.e.\ the region around the zero where $\Psi=0$,
decreases as
$j$ increases. The parameter $\delta$ can be viewed
as a function of the quantization {\em density} (see \cite{elia.mitter.tac01}),
and we do not assume any constraint on its value
(cf.\ the remark below to see why, for open-loop unstable systems, assuming
$\delta\in (0,1)$ does not appear to be restrictive).
The closed-loop system
 is
\be\label{sq}\ba{rcl}
\dot x &=& F(x,\mu)+G(x,\mu)\zeta\\
\dot \zeta &=& q(x,\zeta,\mu)-b(x,\zeta,\mu)\Psi(\bar k \zeta)\;.
\ea\ee
The following quantities are useful below:
\be\label{ub}\ba{l}
\bar w= \dst\max_{(x,\zeta)\in \Omega_{c^2+d^2+1}\,,\,\mu\in \mathcal{P}}
|w(x,\zeta,\mu)|\;,\;
\bar b=\dst\max_{(x,\zeta)\in \Omega_{c^2+d^2+1}\,,\,\mu\in \mathcal{P}} |b(x,\zeta,\mu)|
\;,\; \bar \zeta=\dst\max_{(x,\zeta)\in \Omega_{c^2+d^2+1}} |\zeta|
\;.
\ea\ee
Observe that
the vector field on the right-hand side of
 (\ref{sq}) is discontinuous
and solutions of the system
must be intended in some
generalized sense. In this section,  we focus on Krasowskii solutions, but other
types of solutions are possible (see e.g.\
\cite{ceragioli.depersis.scl07} and references therein).
The main
reason to consider Krasowskii solutions lies in the fact that a
rather complete Lyapunov theory for the study of the stability of
these solutions is available.

\begin{definition}
A curve  $\varphi : [0,+\infty )\to \R ^n$ is a
{\em Krasowskii solution} of
a system of ordinary differential equations
$\dot x=G(t,x)$,
where $G:[0, +\infty ) \times \R ^n \to \R ^n$,
if it is absolutely continuous and for almost every
$t\geq 0$ it satisfies  the  differential inclusion $\dot x\in
K(G(t,x))$, where $K(G(t,x))=\cap _{\delta >0}\overline{\rm
co}\,G(t,B_\delta(x))$ and $\overline{\rm
co}\,G$ is the convex closure of the set $G$. $\triangleleft$
\end{definition}

\noindent In the present case, Krasowskii solutions are absolutely continuous
functions which satisfy the differential inclusion (see e.g.\
\cite{B1},
\cite{ceragioli.depersis.scl07})
\[\ba{l}
\left(\ba{c}
\dot x \\ \dot \zeta
\ea\right)
\in
\left(\ba{c}
F(x,\mu)+G(x,\mu)\zeta \\ q(x,\zeta,\mu)
\ea\right)+
\left\{
\left(\ba{c}
0 \\ -b(x,\zeta,\mu)
\ea\right)v\;,\; v\in K(\Psi(\bar k \zeta))
\right\}
\ea\]
where
\[
\ba{l}
K(\Psi(\bar k \zeta))\subseteq
\left\{\ba{ll}
\{(1+\lambda\delta)\bar k \zeta\;,\;
\lambda\in [-1,1]\} &
\frac{u_j}{1+\delta}<|\bar k \zeta|\le \frac{u_0}{1-\delta}\\[2mm]
\{\lambda(1+\delta)\bar k \zeta\;,\;
\lambda\in [0,1]\} & \frac{u_j}{1+\delta}\ge |\bar k \zeta|\;.
\ea
\right.
\ea
\]

\begin{remark}
Assuming $\delta\in(0,1)$ does not
appear to be restrictive. In fact,
in the differential inclusion above, because of the quantization,
the ``high-frequency" gain becomes
$b(x,\zeta,\mu)(1+\lambda\delta)$, with $\lambda\in [-1,1]$ (similar
arguments hold for the case when the gain is
$b(x,\zeta,\mu)\lambda(1+\delta)$, $\lambda\in [0,1]$). If we allow
$\delta$ to be larger than $1$, then we should design a control
law which stabilizes the system in the presence of
a high-frequency gain which takes any value in a set which includes the zero,
a task which is considerably difficult, if not impossible,  to accomplish
for open-loop unstable systems.  $\triangleleft$
\end{remark}

Then we claim the following version of the so-called ``semi-global backstepping lemma" in
\cite{teel.praly.sicon95} with quantized feedback:
\begin{proposition}\label{l1}
For any $\delta\in (0,1)$, there
exist positive numbers $k^\ast$, $j^\ast$, and $u_0$
such that, for any gain $\bar k\ge k^\ast$ and any number of
quantization levels
$j\ge j^\ast$,  any Krasowskii solution $\varphi$ of the system (\ref{sq}) is
such that, if $\varphi(0)\in \Omega_{c^2+d^2+1}$, then there exists
$T>0$ such that $\varphi(t)\in \Omega_{\sigma}$ for all $t\ge T$.
\end{proposition}
\begin{proof}
Let
\[\ba{l}
k^\ast=\dst\frac{d+1}{d}\dst\frac{\bar w}{b_0}
\dst\frac{1}{\eta (1-\delta)}\;,\;j^\ast =
\left\lceil
\log\left(
\dst\frac{d^2}{(c^2+d^2+d+1)^2}\dst\frac{\eta}{4}\dst\frac{b_0}{\bar b}
\right)
\log\left(
\dst\frac{1-\delta}{1+\delta}
\right)^{-1}
\right\rceil
\ea\]
fix $\bar k\ge k^\ast$ and $j\ge j^\ast$, and choose $u_0=(1+\delta)\bar k
\max_{(x,\zeta)\in \tilde S} |\zeta|$.
To prove convergence of the state to $\Omega_\sigma$,
we need to prove (\cite{ceragioli.depersis.scl07}) that, for any $(x,\zeta)\in S$,
for
any $v\in K(\Psi(\bar k \zeta))$,
\be\label{x1}
\ba{rcl}
&&\dot W(x,\zeta)\\
&=&
\dst\frac{c (c+1)}{(c+1-V(x))^2}
\dst\frac{\partial V}{\partial
x}(F(x,\mu)+G(x,\mu)\zeta)+
\dst\frac{d (d+1)}{(d+1-\zeta^2)^2}2\zeta\cdot
(q(x,\zeta,\mu)-b(x,\zeta,\mu)v)\\
&\le&
\dst\frac{c}{c+1}\dst\frac{\partial V}{\partial
x}F(x,\mu)
-
2\dst\frac{d (d+1)}{(d+1-\zeta^2)^2}\bar k
b(x,\zeta,\mu)\zeta^2+w(x,\zeta,\mu)\zeta-\\[3mm]
&& 2\zeta\dst\frac{d (d+1)}{(d+1-\zeta^2)^2} b(x,\zeta,\mu)[v- \bar k
\zeta]<0\;.
\ea
\ee
Note first that, if $(x,\zeta)\in \tilde S$, then $|\bar k\zeta| \le
u_0$. Hence, depending on the number $j\ge j^\ast$ of quantization
levels, two cases are possible, namely that  the set
$\hat S:=\{(x,\zeta)\in S\,:\, u_j/(1+\delta)<|\bar k\zeta|\le
u_0\}$
is {\em strictly} contained in $\tilde S$ or it is not
.
Consider the former case.
By definition of $K(\Psi(\bar k\zeta))$,
we have $v-\bar k\zeta= \lambda
\delta \bar k\zeta$, and therefore, for $(x,\zeta)\in \hat S\subset\tilde S$,
the inequality above rewrites as (recall (\ref{ub}))
\[
\ba{rcl}
\dot W(x,\zeta)&\le& \dst\frac{c}{c+1}\dst\frac{\partial V}{\partial
x}F(x,\mu)
-2\dst\frac{d (d+1)}{(d+1-\zeta^2)^2}(1+\lambda\delta)\cdot
\bar k
b(x,\zeta,\mu)\zeta^2+w(x,\zeta,\mu)\zeta\\[2mm]
&\le &
\dst\frac{c}{c+1}\dst\frac{\partial V}{\partial
x}F(x,\mu)
-2\dst\frac{d}{d+1}(1-\delta)\bar k
b_0 \zeta^2+ \bar w|\zeta|\;.
\ea
\]
The choice of $k^\ast$ above gives $\dot W(x,\zeta)\le -\bar w \eta$.
Consider now the subset of points in $\tilde S$ such that $|\bar k\zeta| \le
u_j/(1+\delta)$. Such a set is non-void because $\hat S\subset\tilde S$ by hypothesis.
For these points, we have
\[|v-\bar k\zeta|\le |\bar k\zeta|\le
u_j/(1+\delta)\le
\dst\frac{u_0}{1+\delta}\left(
\dst\frac{1-\delta}{1+\delta}
\right)^j
\;,
\]
and the bound on $\dot W(x,\zeta)$ becomes
\[
\dot W(x,\zeta)\le
\dst\frac{c}{c+1}\dst\frac{\partial V}{\partial
x}F(x,\mu)
-\dst\frac{d}{d+1} \bar k
b_0 \zeta^2+
2\dst\frac{u_0}{1+\delta}\left(
\dst\frac{1-\delta}{1+\delta}
\right)^j\dst\frac{(c^2+d^2+d+1)^2}{d(d+1)} \bar b |\zeta|
\;.
\]
The choice of $u_0$ and $j^\ast$ guarantees that for $j\ge j^\ast$
the last term on the right-hand side of the inequality is not larger
than the second term, and this gives
$\dot W(x,\zeta)\le
-\frac{1}{2}\frac{d}{d+1} \bar k
b_0 \eta^2\le - \frac{\bar w\eta}{2(1-\delta)}
$.

Consider now the case when  $\hat S\not\subset \tilde S$, and let
first
$(x,\zeta)\in \hat S\cap\tilde S$. This case is the same as
$(x,\zeta)\in \hat S$ when $\hat S\subset \tilde S$.
Then as before
\[\ba{l}
\dot W(x,\zeta)\le
\dst\frac{c}{c+1}\dst\frac{\partial V}{\partial
x}F(x,\mu)+w(x,\zeta,\mu)\zeta
- 2\dst\frac{d}{d+1}(1-\delta)\bar k
b_0\zeta^2<0\;.
\ea\]
On the other hand, if $(x,\zeta)\not \in \hat S\cap\tilde S$, then necessarily $(x,\zeta)\in
U$. Then we have
\[
\dot W(x,\zeta)\le
\dst\frac{c}{c+1}\dst\frac{\partial V}{\partial
x}F(x,\mu)+w(x,\zeta,\mu)\zeta-
\left\{
\ba{ll}
2\dst\frac{d}{d+1}(1-\delta)\bar k
b(x,\zeta,\mu)\zeta^2, &  \frac{u_j}{1+\delta}<|\bar k\zeta|\le
u_0\\
0\;, & \frac{u_j}{1+\delta}\ge |\bar k\zeta|\;.
\ea\right.
\]
Since the sum of the first two terms on the right-hand side
is negative because $(x,\zeta)\in
U$, and the third term is always
non-positive, we see that $\dot W(x,\zeta)<0$. This ends the proof. 
\end{proof}

\begin{remark}
The constant $k^\ast$ differs from the one in
\cite{teel.praly.sicon95}, \cite{isidori.book.vol2} by the presence of the factor
$(1-\delta)^{-1}$. That is, as expected, the error due to quantization is
counteracted by raising the controller gain. Furthermore, it
is interesting to observe that the constant $j^\ast$, that is the number
of quantization levels,  only depends on
the size of the domain of attraction and of the target set.
\end{remark}

\section{An estimate on the bandwidth}\label{sec.band}

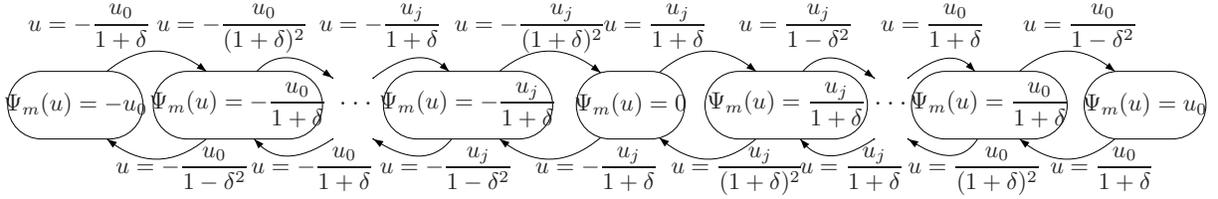
\begin{figure}
\setlength{\unitlength}{1mm}
\begin{center}
\scalebox{0.9}{
\begin{picture}(128,20)(14,-62)


\node[NLangle=0.0,Nw=20.0,Nh=10.0,Nmr=10.0](n0)(2,-48.0){\footnotesize{$\Psi_m(u)=-u_0$}}

\node[NLangle=0.0,Nw=25.0,Nh=10.0,Nmr=10.0](n3)(26.0,-48.0){\footnotesize{$\Psi_m(u)=-\dst\frac{u_0}{1+\delta}$}}

\drawedge[curvedepth=8.0](n0,n3){}

\put(-5,-37){\footnotesize{$u=-\dst\frac{u_0}{1+\delta}$}}

\drawedge[curvedepth=8.0](n3,n0){}

\put(8,-58){\footnotesize{$u=-\dst\frac{u_0}{1-\delta^2}$}}

\node[NLangle=0.0,Nw=0.1,Nh=0.1,Nmr=10.0](n43)(40.0,-44.0){}

\drawedge[curvedepth=5.0](n3,n43){}

\put(14,-37){\footnotesize{$u=-\dst\frac{u_0}{(1+\delta)^2}$}}

\node[NLangle=0.0,Nw=0.1,Nh=0.1,Nmr=10.0](n34)(40.0,-53.0){}

\drawedge[curvedepth=5.0](n34,n3){}

\put(28,-58){\footnotesize{$u=-\dst\frac{u_0}{1+\delta}$}}

\node[NLangle=0.0,Nw=0.1,Nh=0.1,Nmr=10.0](n10)(46.0,-44.0){}

\node[NLangle=0.0,Nw=0.1,Nh=0.1,Nmr=10.0](n11)(46.0,-53.0){}

\put(41,-48){$\ldots$}

\node[NLangle=0.0,Nw=25.0,Nh=10.0,Nmr=10.0](n4)(60.0,-48.0)
{\footnotesize{$\Psi_m(u)=-\dst\frac{u_j}{1+\delta}$}}

\drawedge[curvedepth=5.0](n10,n4){}
\put(38,-37){\footnotesize{$u=-\dst\frac{u_{j}}{1+\delta}$}}

\drawedge[curvedepth=5.0](n4,n11){}
\put(47,-58){\footnotesize{$u=-\dst\frac{u_{j}}{1-\delta^2}$}}

\node[NLangle=0.0,Nw=16.0,Nh=10.0,Nmr=10.0](n5)(84.0,-48.0)
{\footnotesize{$\Psi_m(u)=0$}}

\drawedge[curvedepth=8.0](n4,n5){}

\put(58,-37){\footnotesize{$u=-\dst\frac{u_j}{(1+\delta)^2}$}}

\drawedge[curvedepth=8.0](n5,n4){}

\put(70,-58){\footnotesize{$u=-\dst\frac{u_j}{1+\delta}$}}

\node[NLangle=0.0,Nw=24.0,Nh=10.0,Nmr=10.0](n6)(107.0,-48.0)
{\footnotesize{$\Psi_m(u)=\dst\frac{u_j}{1+\delta}$}}

\node[NLangle=0.0,Nw=0.1,Nh=0.1,Nmr=10.0](n65)(120.0,-44.0){}
\drawedge[curvedepth=5.0](n6,n65){}
\put(100,-37){\footnotesize{$u=\dst\frac{u_j}{1-\delta^2}$}}

\node[NLangle=0.0,Nw=0.1,Nh=0.1,Nmr=10.0](n56)(120.0,-53.0){}
\drawedge[curvedepth=5.0](n56,n6){}
\put(109,-58){\footnotesize{$u=\dst\frac{u_j}{1+\delta}$}}

\drawedge[curvedepth=8.0](n5,n6){}

\put(80,-37){\footnotesize{$u=\dst\frac{u_j}{1+\delta}$}}

\drawedge[curvedepth=8.0](n6,n5){}

\put(90,-58){\footnotesize{$u=\dst\frac{u_j}{(1+\delta)^2}$}}

\put(120,-48){$\ldots$}

\node[NLangle=0.0,Nw=23.0,Nh=10.0,Nmr=10.0](n7)(137.0,-48.0)
{\footnotesize{$\Psi_m(u)=\dst\frac{u_0}{1+\delta}$}}

\node[NLangle=0.0,Nw=0.1,Nh=0.1,Nmr=10.0](n75)(125.0,-44.0){}
\drawedge[curvedepth=5.0](n75,n7){}
\put(121,-37){\footnotesize{$u=\dst\frac{u_0}{1+\delta}$}}

\node[NLangle=0.0,Nw=0.1,Nh=0.1,Nmr=10.0](n57)(125.0,-53.0){}
\drawedge[curvedepth=5.0](n7,n57){}
\put(125,-58){\footnotesize{$u=\dst\frac{u_0}{(1+\delta)^2}$}}

\node[NLangle=0.0,Nw=18.0,Nh=10.0,Nmr=10.0](n8)(160.0,-48.0)
{\footnotesize{$\Psi_m(u)=u_0$}}

\drawedge[curvedepth=8.0](n7,n8){}

\put(142,-37){\footnotesize{$u=\dst\frac{u_0}{1-\delta^2}$}}

\drawedge[curvedepth=8.0](n8,n7){}

\put(146,-58){\footnotesize{$u=\dst\frac{u_0}{1+\delta}$}}

\end{picture}
}
\caption{\label{fig.multivalued2}
The graph
illustrates how the function
$\Psi_{m}(u)$ takes values depending on $u$, $u=\bar k \zeta$. Each edge
connects two nodes, and is
labeled with the condition ({\it guard}) which triggers the transition
from the starting node to the destination
node.
}
\end{center}
\end{figure}

The stabilization technique examined in the previous section has two main ingredients:
the selection of the set of control values, and the switching law
which schedules them. As a possible performance measure
of the control law can then be taken the number of times the controller
switches to a new value within an interval of time divided by the length of
the interval itself. Such a measure is given the name of
{\em bandwidth}, because of the obvious implication in the scenario
in which the quantized control is fed back to the process through a
finite bandwidth communication channel.
The Krasowskii solutions considered above do not exclude the
possibility to have accumulation of switching points in finite time.
To circumvent the possible occurrence of chattering or Zeno phenomenon,
we introduce a modified quantizer following
\cite{hayakawa.et.al.acc06}. The modified quantizer is obtained from
(\ref{ppsi}) to which,  for each quantization level, a new
one is added, to obtain the following multi-valued map
(see \cite{hayakawa.et.al.acc06}, Figure 3.4, for a pictorial
representation of the map):
\be\label{ppsi.haya}
\Psi_m(u)=
\left\{\ba{ccclrcll}
u_i &\dst\frac{1}{1+\delta}u_i& < & u &\le &
\dst\frac{1}{1-\delta}u_i\;, &  0\le i\le j\\[3mm]
\dst\frac{1}{1+\delta}u_i & \dst\frac{1}{(1+\delta)^2}u_i& < & u & \le&
\dst\frac{1}{(1+\delta)(1-\delta)}u_i\;, & 0\le i\le j\\[3mm]
0 & 0&\le & u &\le & \dst\frac{1}{1+\delta}u_j & \\[3mm]
-\Psi_m(-u) & & & u & < & \hspace{-1cm} 0\;. &
\ea\right.
\ee
Since the map above is multi-valued, we need to specify the law
according to which $\Psi_m(u(t))$ changes its value as $u(t)$ evolves.
This law is illustrated by the graph in Fig.\
\ref{fig.multivalued2}.
At time $t=0$, we set $\Psi_m(u(0))=\Psi(u(0))$. This value of
$\Psi_m(u(0))$ identifies a node of the graph. If the value of $u(0)$
fulfills one of the conditions of the edges
leaving the node, then a transition is triggered and
the quantizer takes the new value -- which is denoted by
$\Psi_m(u(0^+))$ --  given by the destination node.
For $t>0$, $\Psi_m(u(t))$  remains constant until
$u(t)$ triggers a transition of $\Psi_m(u(t))$
to the new value, denoted by $\Psi_m(u(t^+))$, again chosen according to the graph
of Fig.\ \ref{fig.multivalued2}.
We refer to \cite{hayakawa.et.al.acc06}, Section 3,  for further
details on the switching mechanism.

The first
observation is that the result proven in the previous section
continues to hold even in the presence of the modified quantizer. As
a matter of fact, Proposition \ref{l1} was proven showing that the derivative
\be\label{derivative}
\left(\ba{cc}
\dst\frac{\partial W}{\partial x} & \dst\frac{\partial W}{\partial \xi}
\ea\right)
\left(\ba{c}
F(x,\mu)+G(x,\mu)\zeta \\ q(x,\zeta,\mu)
-b(x,\zeta,\mu)v
\ea\right)
\ee
was strictly negative for all $(x,\zeta)\in S$, $\mu\in \mathcal{P}$
and  $v\in K(\Psi(\bar k\zeta))$. Now, if we adopt the
modified quantizer defined above, the closed-loop  system becomes
the switched system
\be\label{sqm}\ba{rcl}
\dot x &=& F(x,\mu)+G(x,\mu)\zeta\\
\dot \zeta &=& q(x,\zeta,\mu)-b(x,\zeta,\mu)\Psi_m(\bar k \zeta)\;,
\ea\ee
and proving stability of the  (unique)
solution amounts simply to show that (\ref{derivative}) is still
negative
when $v$ is replaced by $\Psi_m(\bar k \zeta)$. This is an immediate
consequence of the result in the previous section:
\begin{corollary}\label{c1}
For any $\delta\in (0,1)$, there
exist positive numbers $k^\ast$, $j^\ast$, and $u_0$
such that, for any gain $\bar k\ge k^\ast$ and any number of
quantization levels
$j\ge j^\ast$, the unique solution $\varphi$ of
of the system (\ref{sqm}),
is
such that, if $\varphi(0)\in \Omega_{c^2+d^2+1}$, then there exists
$T>0$ such that $\varphi(t)\in \Omega_{\sigma}$ for all $t\ge T$.
\end{corollary}

\begin{proof}
By definition of
$\Psi_m$ and $K(\Psi(\bar k\zeta))$, for each $|\bar k \zeta|\le u_0(1-\delta)^{-1}$, $\Psi_m(\bar k
\zeta)\in K(\Psi(\bar k\zeta))$. This ends the proof.
\end{proof}

Now we make the notion of bandwidth more precise.
Let first $0=t_0,t_1,t_2,\ldots$ be
the sequence of switching times, that is the times at which
the control law $u=-\Psi_m(\bar k \zeta)$ changes its value,  and
$B(t_\kappa)$ the number of quantization levels used to encode the
control at time $t_\kappa$
(but we could equivalently use the number of bits
used to encode the value transmitted at time $t_k$).
For each
$t$, let $\kappa$ the largest integer for which $t\ge t_\kappa$,
and define the average data rate over the
interval $[t_0, t]$ as the quantity
$R_{av}[t_0, t]=(\sum_{j=0}^\kappa B(t_j))/(t-t_0)$,
that is the total number of values taken by the quantized control on  the interval
$[t_0,t]$, divided by the duration of the interval. Moreover,
we denote as the average data rate the limit
$R_{av}={\lim\sup}_{t\to \infty} R_{av}[t_0, t]={\lim\sup}_{t\to \infty}
(\sum_{j=0}^\kappa B(t_j))/(t-t_0)$.
Under the conditions of Corollary \ref{c1}, the following result provides
a bound on the average data rate needed to guarantee semi-global
practical stabilization. In the statement, we refer to the following
quantities
\be\label{notation.bw}
\bar q =
\dst\max_{(x,\zeta)\in \Omega_{c^2+d^2+1},\mu\in\mathcal{
P}}|q(x,\zeta,\mu)|
\;,\;
k_0 = \dst\frac{d+1}{d}\dst\frac{\bar w}{b_0}
\dst\frac{1}{\eta}\;.
\ee
\begin{proposition}
Let $\bar k=k^\ast$. The unique solution $\varphi$ of
the system (\ref{sqm})
is
such that, if $\varphi(0)\in \Omega_{c^2+d^2+1}$, then there exists
$T>0$ such that $\varphi(t)\in \Omega_{\sigma}$ for all $t\ge T$,
and the associated average data rate is not greater than
$
\frac{4j+1}{DT_m}
$,
where
$
DT_m =\frac{1}
{(\bar \zeta\rho^{j-1})^{-1}\bar q+k_0\bar b}
\frac{\delta}{1+\delta}
$.
\end{proposition}

\begin{proof}
The proof boils down to estimate a lower bound on the
time which elapses between two consecutive switching times
({\em inter-switching} time). The estimate is found by focusing on
the function of time $\Psi_m(k\zeta(t))$, and in particular on the
value $|\Psi_m(k\zeta(t^+))|$ and where $t^+$
denotes $|\Psi_m(k\zeta(t))|$ soon after the
switching. Then, the smallest distance to be covered by $k\zeta(t)$
before a new switching takes place, and the largest velocity at
which the function $k\zeta(t)$ evolves are computed. It is enough to carry out
these calculations in three cases.
If
$|\Psi_m(k\zeta(t^+))|=u_i$, with $0\le i\le j$, then $|u|$ remains
equal to $u_i$ for not less than
\[
\dst\frac{\dst\frac{u_i}{1+\delta}\dst\frac{\delta}{1-\delta}}
{\bar k(\bar q+ \bar b u_i)}
= \dst\frac{{\bar \zeta}}
{\bar q+\bar k(1-\delta)\bar b
\bar \zeta\rho^{i-1}}
\rho^{i}\dst\frac{\delta}{1-\delta}\;,
\]
where we have exploited the definition of $u_i$, $u_0$, $\rho=(1-\delta)(1+\delta)^{-1}$,
and (\ref{notation.bw}).
Observe
that by the definition of $k^\ast$ in  Proposition \ref{l1},
$\bar k(1-\delta)=k^\ast(1-\delta)= k_0$.
Hence,  the bound above becomes equal to
\be\label{ub2}\ba{c}
\dst\frac{{\bar \zeta}}
{\bar q+k_0\bar b
\bar \zeta\rho^{i-1}}
\rho^{i}\dst\frac{\delta}{1-\delta}
\;.
\ea\ee
With similar arguments it can be shown that a lower bound on the time $|u|$ remains
equal to $u_i(1+\delta)^{-1}$ is given by:
$
\frac{{\bar \zeta}}
{\bar q+k_0\bar b
\bar \zeta\rho^{i-1}\dst\frac{1}{1+\delta}}
\rho^{i}\frac{\delta}{1-\delta}
$.
Finally a lower bound on the dwell time when $|\Psi_m(k\zeta(t^+))|=0$
is
$
\frac{\bar \zeta}
{\bar q}\rho^{j}\dst\frac{2}{1+\delta}
$.
Hence, comparing the three estimates above it is seen
that the inter-switching time
can not be less than (\ref{ub2}) with $i=j$. Indeed, the bound (\ref{ub2})
is a monotonically
decreasing function of $i$ and hence the minimum is reached at $i=j$ and is
equal to
\[
\ba{l}
DT_m=
\dst\frac{{\bar \zeta}}
{\bar q+k_0\bar b
\bar \zeta\rho^{j-1}}
\rho^{j}\dst\frac{\delta}{1-\delta}=
\dst\frac{1}
{\dst\frac{\bar q}{\bar \zeta\rho^{j-1}}+k_0\bar b}
\dst\frac{\delta}{1+\delta}
\;.
\ea\]
Let now $t\in [t_\kappa,t_{\kappa+1})$, $\kappa\ge 1$. Then $\kappa+1$, the number of switchings
in the interval $[t_0,t]$, satisfies
$
\kappa\le \dst\frac{t-t_0}{DT_m}
$. Then
\[\ba{rcl}
R_{av}&=&{\lim\sup}_{t\to \infty} R_{av}[t_0, t]={\lim\sup}_{t\to \infty}
\dst\frac{\sum_{\ell=0}^\kappa (4j+1)}{t-t_0}={\lim\sup}_{t\to \infty}
\dst\frac{(\kappa+1)(4j+1)}{t-t_0}\\
&=&{\lim\sup}_{t\to \infty}
(\dst\frac{\kappa}{t-t_0}+\dst\frac{1}{t-t_0})(4j+1)\le{\lim\sup}_{t\to \infty}
(\dst\frac{1}{DT_m}+\dst\frac{1}{t-t_0})(4j+1)
=\dst\frac{4j+1}{DT_m}
\;.
\ea\]
\end{proof}

\begin{remark}
A special case is when $q(x,\zeta,\mu)$ is identically zero. In this
case it is seen that $|\bar k\zeta(t)|$ is a monotonically decreasing
function as far as $u_j(1+\delta)^{-2}<|\bar k\zeta(t)|\le
u_0(1-\delta)^{-1}$. Hence, at some finite time $\bar t$, the control becomes equal
to zero, and $\zeta(t)$ remains equal to the value $\zeta(\bar t)$
for all $t\ge \bar t$. As a result, switching stops, and for $t\ge \bar t$, the number
$\kappa$ of switchings in any
interval of time $[t,t_0]$ remains constant and finite.
Hence, $R_{av}=0$. We conclude that in the special case
$q(x,\zeta,\mu)=0$, the control law $u=-\Psi_m(\bar k\zeta)$ guarantees
semi-global practical stability with an average data rate equal to
zero. $\triangleleft$
\end{remark}

\section{Ternary controller}\label{sec.ternary}

\begin{figure}
\begin{center}
\begin{picture}(78,25)(0,-22)

\put(26,1.5){\makebox(0,0)[cc]{$\zeta\le \eta/2$}}
\put(54,1.5){\makebox(0,0)[cc]{$\zeta\le -\eta$}}
\put(54,-21){\makebox(0,0)[cc]{$\zeta\ge -\eta/2$}}
\put(26,-21){\makebox(0,0)[cc]{$\zeta\ge \eta$}}

\node[NLangle=0.0,Nw=14,Nh=14,Nmr=6.22](n0)(68.0,-10.0){$u=\bar k$}

\node[NLangle=0.0,Nw=14,Nh=14,Nmr=6.22](n1)(40.0,-10.0){$u=0$}

\node[NLangle=0.0,Nw=14,Nh=14,Nmr=6.22](n2)(12.0,-10.0){$u=-\bar k$}

\drawedge[curvedepth=8.0](n2,n1){}

\drawedge[curvedepth=8.0](n1,n0){}

\drawedge[curvedepth=8](n0,n1){ }

\drawedge[curvedepth=8](n1,n2){}

\end{picture}
\end{center}
\caption{\label{fig.graph} The ternary switched controller.}
\end{figure}
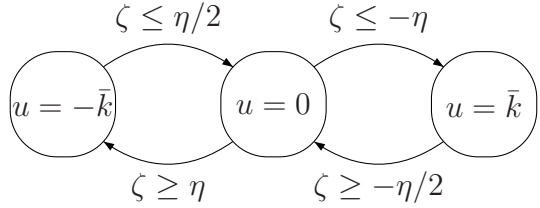

In this section we remark that Proposition \ref{l1} can be also obtained
using a ternary controller. Let $\eta$ be the  positive constant introduced
at the end of Section \ref{sec.prel}
in the definition of the neighborhood $U$,
and introduce the following sets:
\[\ba{l}
\Omega_-=\{
(x,\zeta)\in \Omega_{c^2+d^2+1}\,:\, \zeta\ge \eta
\}\;,\;
\Omega_0=\{
(x,\zeta)\in \Omega_{c^2+d^2+1}\,:\, |\zeta|< \eta
\}\;,\\
\Omega_+=\{
(x,\zeta)\in \Omega_{c^2+d^2+1}\,:\, \zeta\le -\eta
\}\;.
\ea\]
These sets are depicted in Figure \ref{fig.sets}. Assume without loss of
generality that $\eta$ is small enough that $\Omega_-,\Omega_+$ are not
void. We propose the following controller. (Similar elementary
 controllers have been studied in
\cite{kaliora.astolfi.tac04} for a different class of nonlinear
systems.)
At the initial time $t=0$, assume that $(x,\zeta)\in
\Omega_{c^2+d^2+1}$, and set the control value as
\be\label{control.0}
u(0)=\left\{\ba{rcl}
-\bar k & if & \zeta(0)\ge \eta\\
0 & if & |\zeta(0)|< \eta\\
\bar k & if & \zeta(0)\le -\eta\;.
\ea\right.
\ee
As in the previous section, for $t\ge 0$,
the controller is chosen
according to the law
\be\label{control}
u(t^+)=\left\{\ba{rcl}
-\bar k & if & ([u(t)=0]\wedge [\zeta(t)\ge \eta])\vee
([u(t)=-\bar k]\wedge [\zeta(t)> \eta/2])\\
0 & if &
\ba{l}
([u(t)=-\bar k]\wedge [\zeta(t)\le \eta/2])\vee
([u(t)=\bar k]\wedge [\zeta(t)\ge
-\eta/2])\\~
\vee ([u(t)=0]\wedge [|\zeta(t)|<
\eta])
\ea
\\
\bar k & if & ([u(t)=0]\wedge [\zeta(t)\le -\eta])
\vee ([u(t)=\bar k]\wedge [\zeta(t)< -\eta/2])
\;,
\ea\right.
\ee
with $\bar k>0$ a parameter to design. This law could also be
described by an automaton analogous to the one in Fig.\ \ref{fig.multivalued2}
but with three states only (see Fig.\ \ref{fig.graph}).

The stability result with the ternary controller reads as follows:
\begin{proposition}\label{p1}
There exists a choice of $\bar k$ such that
the Lyapunov function $W(x,\zeta)$, computed along any
trajectory of the closed-loop system
(\ref{basic}), (\ref{control.0}), (\ref{control})
which starts in  $S$,
satisfies
$\dot W(x(t),\zeta(t))<0$
for all $(x(t),\zeta(t))\in S$.
\end{proposition}

\begin{proof}
It has already proven that, for $(x,\zeta)\in U$, $\dot W(x(t),\zeta(t))<0$
with $u=0$.
On the other hand, for $u=-\bar k{\rm sgn}\,\zeta$,
\[\ba{rcl}
\dot W(x,\zeta,u)&=&
\dst\frac{c (c+1)}{(c+1-V(z))^2}\dst\frac{\partial V}{\partial
x}F(x,\mu)+
w(x,\zeta,\mu)\zeta-2\dst\frac{d (d+1)}{(d+1-\zeta^2)^2}\zeta
b(x,\zeta,\mu)\bar k|\zeta|\\
&\le&\dst\frac{c}{c+1}\dst\frac{\partial V}{\partial
x}F(x,\mu)+ |w(x,\zeta,\mu)|\,|\zeta|-2\bar k b_0 \dst\frac{d}{d+1}|\zeta|
\;,
\ea
\]
and choosing $\bar k \ge \frac{1}{b_0}\frac{d+1}{d} \bar w$,
with $\bar w$ as in (\ref{ub}),
we have, for all $(x,\zeta)\in \Omega_{c^2+d^2+1}$ such that $|\zeta| \ge \eta$,
$
\dot W(x,\zeta,u)\le -\bar w |\zeta|\le - \bar w\eta
$.
This concludes the proof.
\end{proof}

\begin{remark}
As for quantized control, it is possible to give an estimate on the
bandwidth needed to implement the ternary controller. Using the same
arguments of Section \ref{sec.band}, one can show that an upper bound
on the average data rate is
$6(\bar q\eta^{-1}+\bar b k_0)$
where $k_0$ is the quantity defined in (\ref{notation.bw}).
The two estimates on the data rate using quantized and ternary control
are similar but it is difficult to compare the two control laws.
Arguably, in some cases, the ternary control above
is easier to implement than the quantized control (cf.\ \cite{li.baillieul.tac04}
to see how binary control is more robust to
changes  for linear scalar systems).
On the other hand, while the state is approaching the target set,
the amplitude of the changes in the values taken by the
quantized controller are less pronounced  than
the ternary control and in some cases, to keep the state in the
vicinity of the origin, the quantized controller may
require a minor effort than the ternary controller.
\end{remark}

\section{Applications}\label{sec.appl}
In this section we emphasize a number of cases to which the
previous results apply.

\subsection{Systems with uniform relative degree $\ge 1$}\label{subs.urd}

It is well-known that a nonlinear input-affine system is said to have a {\em
uniform} relative degree $r$ if it has a relative degree $r$ at
$x_0$ for each $x_0\in \mathbb{R}^n$. It is also well-known that
there exists a globally defined diffeomorphism which changes the
system into one of the following form (see e.g.\ Proposition 9.1.1.\
in \cite{isidori.book.vol1}):
\be\label{s}\ba{rcl}
\dot z &=& f(z,\xi_1)\\
\dot \xi_i &=& \xi_{i+1} \;,\; 1\le i\le r-1\\
\dot \xi_{r} &=& \bar q(z,\xi)+\bar b(z,\xi)u
\ea\ee
with $z\in \R^{n-r}$, and $\bar b(z,\xi)\ge b_0>0$ for all $(z,\xi)$.
Systems
like the one above restricted to the components
$z,\xi_1,\ldots,\xi_{r-1}$, with $\xi_r$ viewed as an input, can be always
stabilized by means of a linear high-gain partial-state feedback
(\cite{isidori.book.vol1}, Theorem 9.3.1), provided that the origin $z=0$
is a  globally asymptotically stable equilibrium point for $\dot z=f(z,0)$, i.e.
system (\ref{s}) is minimum-phase. As a matter of fact,
for any $R>0$, there exists a linear ``control law" $\xi_r=-a\xi$,
with $a$ a row vector depending on $R$ and
$\xi=(\xi_1\ldots\xi_{r-1})^T$, such that
every solution of
\be\label{ss}\ba{rcl}
\dot z &=& f(z,\xi_1)\\
\dot \xi_i &=& \xi_{i+1} \;,\; 1\le i\le r-2\\
\dot \xi_{r-1} &=& -a\xi
\ea\ee
starting from the
cube in $\R^{n-1}$ whose edges are $2R$ long, asymptotically converges to the
origin. Perform the change of coordinates
$\xi_r=-a\xi+\zeta$,
let $x=(z^T\,\xi^T)^T$,
and rewrite (\ref{s}) as
\[
\ba{rcl}
\dot x &=& F(x)+G\zeta\\
\dot \zeta &=& q(x,\zeta)+b(x,\zeta)u\;,
\ea
\]
%
%
where $F(x)$ is the vector field on the right-hand side of
(\ref{ss}), and $G$, $q$, $b$ are understood from the context. The
system $\dot x=F(x)$ satisfies the ULP property. We conclude that
both Proposition \ref{l1}, Corollary \ref{c1} and Proposition \ref{p1} can be applied to
system (\ref{s}) to obtain
\begin{proposition}\label{pr}
Consider a minimum-phase nonlinear system  of the form (\ref{s}).
For any $R>0$ and any $\varepsilon>0$, there exist quantized
feedback laws $u=-\Psi(\bar k\zeta)$, $u=-\Psi_m(\bar k\zeta)$, or a ternary feedback law
(\ref{control.0}), (\ref{control}), with $\zeta=a\xi+\xi_r$,
and a time $T>0$, such that any trajectory $\varphi$
of the closed-loop system
which starts in the cube centered at the origin of side $2R$
lies in the cube centered at the origin of side
$2\varepsilon$ for all $t\ge T$.
\end{proposition}

\begin{remark}
The Proposition shows that it is as simple as in the non-quantized case
to stabilize
nonlinear minimum-phase systems with {\em quantized measurements}
provided that the {\em relative degree} of the system is {\em one}.
In fact, if this is the case, then $\zeta$ coincides with the output
of the system.
\end{remark}

In the remaining subsections, we shall refer to systems
for which similar results apply
as {\em semi-globally practically stabilizable} systems.

\vspace{-0.1cm}

\subsection{Robust quantized stabilization of nonlinear systems}

\vspace{-0.1cm}

In this section we propose a quantized controller to stabilize
nonlinear systems of the form
\be\label{fff}\ba{l}
\dot z = F(\mu)z+G(\xi_1,\mu) \xi_1\\
\dot \xi_i = q_{i0}(\xi_1,\ldots,\xi_i,\mu)z
+\dst\sum_{j=1}^{r-1}q_{i,j}(\xi_1,\ldots, \xi_{i},\mu)\xi_j
+b_{i}(\xi_1,\ldots,
\xi_{i},\mu)\xi_{i+1}\;,\; 1\le i\le r-1\\
\dot \xi_{r} =
q_{r0}(\xi_1,\ldots, \xi_{r},\mu)z+
\dst\sum_{i=1}^{r}q_{ri}(\xi_1,\ldots, \xi_{r},\mu)\xi_i
+b_{r}(\xi_1,\ldots,
\xi_{r},\mu)u\;,
\ea\ee
where $z\in \mathbb{R}^{n-r}$, and $b_i(z,\xi_1,\ldots,\xi_i,\mu)\ge b_{i0}>0$ for all
$(z,\xi_1,\ldots,\xi_i)\in \mathbb{R}^{n-r+i}$ and $\mu\in {\cal P}$. We also assume
that, for all $\mu\in \mathcal{P}$, there exists $P(\mu)=P^T(\mu)>0$ such that
$F^T(\mu)P(\mu)+P(\mu)F(\mu)\le -I$.
The first fact we recall is the following (\cite{freeman.kokotovic.aut93},
\cite{isidori.book.vol2}):
\begin{lemma}\label{lemma3}
There exists an $(r-1)\times (r-1)$ matrix $M(\xi)$
and a $1\times (r-1)$ vector $\delta(\xi)$ of smooth functions such
that $\xi^T M(\xi)\xi$ is a definite positive and proper function,
and the function
$
V(z,\xi)=z^TP(\mu)z+ \xi^T M(\xi)\xi
$
satisfies $\dot V(z,\xi)|_{(\ref{fff})-u=\delta(\xi)\xi}\le -\varepsilon
V(z,\xi)$, where $\dot V(z,\xi)|_{(\ref{fff})-u=\delta(\xi)\xi}$ denotes the derivative of
(\ref{fff}) in closed loop with $u=\delta(\xi)\xi$.
\end{lemma}
By the change of coordinates $\zeta=\xi_r-\delta(\xi)\xi$, letting
as before $x^T=(z^T,\xi^T)^T$,
it is immediate to see that we are in the setting of
Proposition \ref{l1} or Proposition \ref{p1}, and systems of the form (\ref{fff}) can
be semi-globally practically stabilized by a quantized  or ternary
controller.

\subsection{A simple output-feedback switched stabilization scheme}
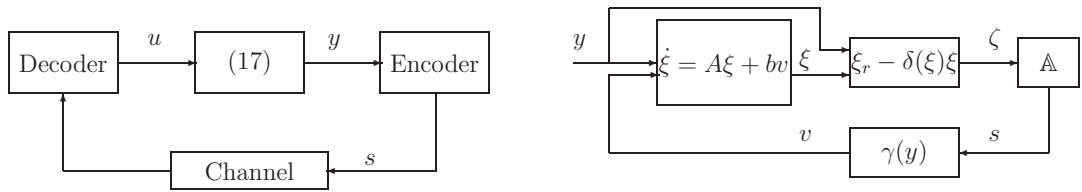
\begin{figure}
\begin{center}
\begin{tabular}{cc}
\scalebox{0.8}{
\begin{picture}(88.75,30)(0,0)
\put(36,4){\framebox(26,6)[cc]{{\rm Channel}}}
\put(40,20){\framebox(18,10)[cc]{(\ref{ts})}}
\put(32,28){$u$}
\put(62,28){$y$}
\put(68,8){$s$}
\put(70.75,20){\framebox(18,10)[cc]{Encoder}}
\put(9,20){\framebox(18,10)[cc]{Decoder}}
\put(27,25){\vector(1,0){13.25}}
\put(58,25){\vector(1,0){13}}
\put(80,20){\line(0,-1){13}}
\put(80,7){\vector(-1,0){18.25}}
\put(36,7){\line(-1,0){18}}
\put(18,7){\vector(0,1){13}}
\end{picture}
}
&
\scalebox{0.8}{
\begin{picture}(90,34)(0,0)
\put(6,28){$y$}
\put(43.5,24.5){$\xi$}
\put(43.5,12){$v$}
\put(75,28){$\zeta$}
\put(75,12){$s$}
\put(20,18){\framebox(22,14)[cc]{$\dot \xi=A\xi+bv$}}
\put(52,21){\framebox(18,8)[cc]{$\xi_r-\delta(\xi)\xi$}}
\put(52,6){\framebox(18,8)[cc]{$\gamma(y)$}}
\put(80,21){\framebox(10,8)[cc]{$\mathbb{A}$}}
\put(6,25){\vector(1,0){14}}
\put(42,23){\vector(1,0){10}}
\put(12,25){\line(0,1){9}}
\put(12,34){\line(1,0){35}}
\put(47,34){\line(0,-1){7}}
\put(47,27){\vector(1,0){5}}
\put(70,25){\vector(1,0){10}}
\put(85,21){\line(0,-1){11}}
\put(85,10){\vector(-1,0){15}}
\put(52,10){\line(-1,0){40}}
\put(12,10){\line(0,1){13}}
\put(12,23){\vector(1,0){8.25}}
\end{picture}
}\end{tabular}
\caption{\label{fx} In the picture on the left,
the switched output feedback controller for system (\ref{ts})
is implemented through a network. The encoder is depicted in the picture on the right.
The block labeled with $\mathbb{A}$ is
the automaton depicted in Fig.\ \ref{fig.multivalued2} (quantized control)
or the automaton described by (\ref{control.0})-(\ref{control}) (ternary control).
The device which converts
the values generated by $\mathbb{A}$ into packets of bits
which can be transmitted through the network
is not depicted for the sake of simplicity. }
\end{center}
\end{figure}

Consider the nonlinear system
\be\label{ts}
\ba{rcl}
\dot x &=& F(\mu)x +G(y,\mu)y+\bar g(\mu)\gamma(y) u\\
\dot y &=& H(\mu)x+K(y,\mu)y\;,
\ea
\ee
with $x\in \mathbb{R}^n$, $y\in \mathbb{R}$ the measured output,
and $\gamma(y)$ a smooth function bounded away from zero.
Under appropriate conditions, namely  (\cite{marino.tomei.I.tac93},
\cite{marino.tomei.II.tac93}, and
also \cite{isidori.book.vol2}, Section
11.3)
(i)  the system has a well-defined uniform
relative degree $r\ge 2$ and (ii) its zero dynamics is globally
asymptotically stable, one can prove that, for the system above, to which it is
appended the additional dynamics
\be\label{ad}
\ba{rcl}
\dot \xi_i &=& -\lambda_{i-1} \xi_i +\xi_{i+1}\;,\; 2\le i\le r-1\\
\dot \xi_r &=& -\lambda_{r-1} \xi_r +\gamma(y)u\;,
\ea
\ee
there exists a change of coordinates $z=T(x,y,\xi,\mu)$,
linear in $(x,y,\xi,\mu)$, which transforms the extended system into
\[
\ba{rcl}
\dot z &=& \tilde F(\mu)z +\tilde G(y,\mu)y\\
\dot y &=& \tilde H(\mu)z+\tilde K(y,\mu)y+b(\mu)\xi_2\\
\dot \xi_i &=& -\lambda_{i-1} \xi_i +\xi_{i+1}\;,\; 2\le i\le r-1\\
\dot \xi_r &=& -\lambda_{r-1} \xi_r +\gamma(y)u\;,
\ea
\]
with $b(\mu)$ bounded away from zero. This system is in the form
(\ref{fff}), and therefore there exists a quantized or a
ternary
controller depending on $y,\xi_2,\ldots, \xi_r$ for it.
The appended dynamics (\ref{ad}) with $u$ given by
(\ref{ppsi}), (\ref{ppsi.haya}) or (\ref{control.0}), (\ref{control}), and
$\zeta=\xi_r-\delta(\xi)\xi$, $\xi=(y,\xi_2, \ldots, \xi_{r-1})$,
is a dynamic output feedback controller which semi-globally
practically stabilizes the system (\ref{ts}).
The implementation of the closed-loop system through a network
in the case of quantized or ternary controller is illustrated in Fig.\ \ref{fx}.
The decoder, on the other hand,
is
a device which carries out the inverse operation with
respect to the encoder, and is not depicted for the sake of simplicity.\\
Compared with \cite{depersis.ijrnc06}, the solution proposed here does not require a
copy of the system to control, and in fact applies to a class
of systems which, although less general than the class in
\cite{depersis.ijrnc06}, present model uncertainty. Moreover, the
dynamics of the ``encoder" on the sensor side is {\em linear},
and therefore it requires less
computational effort than in \cite{depersis.ijrnc06}.
A similar class
as (\ref{ts}) was considered in \cite{liberzon.med07}, where
the output is quantized with no pre-processing. However, in that paper, the
control law $u$ must be designed so as to guarantee input-to-state
stability with respect to state measurements errors, a task which may be
considerably harder than designing the control law as in Lemma
\ref{lemma3}. Observe that we do {\em not} employ a dense
quantization, that is we do not require a small quantization error (the quantization density can
be any number in $(0,1)$) to compensate for the lack of
input-to-state stability.

\section{Conclusion}
We have discussed
results on the problem of stabilizing nonlinear systems using a
finite number of control values and in the presence of parametric
uncertainty. These results are instrumental to solve important
control problems by quantized feedback, a few of which have been presented in the paper.
The tools presented in the paper are suitable to tackle other
important control problems by quantized feedback such as the output
regulation problem, a topic on which future research could focus.


\begin{thebibliography}{99}

\bibitem
{bacciotti.tac89} A.\ Bacciotti.
\newblock
Further remarks on potentially global stabilizability.
\newblock  {\em IEEE Transactions on Automatic Control},  34 (6), 637--639, 1989.

\bibitem
{B1}   Bacciotti  A. Stabilization by means of state space
depending switching rules.  {\em Systems and Control Letters}, {\bf 53},
195-201,  2004.

\bibitem
{ceragioli.depersis.scl07} F.\ Ceragioli and C.\ De Persis.
\newblock
Discontinuous stabilization of nonlinear systems: Quantized and switching controls.
\newblock  {\em Systems \& Control Letters}, Vol. 56, 7-8, 461-473, 2007.


\bibitem
{depersis.ijrnc06}  C.\ De Persis.
\newblock
On stabilization of nonlinear systems under data rate constraints using output measurements.
\newblock  {\em International Journal of Robust and Nonlinear Control}, 16, 315-332, 2006.

\bibitem
{freeman.kokotovic.aut93}
R.A.\ Freeman and P.\ Kokotovic.
\newblock Design of `softer' robust nonlinear control laws.
\newblock {\em Automatica}, 32, 733-746, 1993.

\bibitem
{hayakawa.et.al.acc06} T.\ Hayakawa, H.\ Ishii, and K.\ Tsumura.
\newblock  Adaptive quantized control for nonlinear uncertain systems.
\newblock  In {\em Proc.\ 2006 American Control Conference}, Minneapolis, Minnesota, 2006.

\bibitem
{isidori.book.vol1} A.\ Isidori.
\newblock \emph{Nonlinear Control
Systems}, 3rd edition.
\newblock Springer-Verlag, New York, 1995.


\bibitem
{isidori.book.vol2} A.\ Isidori.
\newblock \emph{Nonlinear Control
Systems}, Vol.\ 2.
\newblock Springer, London, 1999.

\bibitem
{kaliora.astolfi.tac04}
G.\ Kaliora and A.\ Astolfi.
\newblock Nonlinear control of feedforward systems with bounded
signals.
\newblock {\em IEEE Trans.\ Auto.\ Contr.}, 49 (11), 1975--1990, 2004.

\bibitem
{li.baillieul.tac04} K.\ Li and J.\ Baillieul.
\newblock Robust Quantization for Digital Finite Communication Bandwidth (DFCB)
Control.
\newblock {\em IEEE Trans.\ Auto.\ Contr.}, 49, 1573--1584,
2004.

\bibitem
{liberzon.aut03}
D.\ Liberzon.
\newblock Hybrid feedback stabilization of systems with quantized
signals.
\newblock {\em Automatica}, 39: 1543--1554, 2003.

\bibitem
{liberzon.med07}
D.\ Liberzon.
\newblock Observer-based quantized output feedback control of nonlinear systems.
\newblock In {\em Proc.\ MED'07 }, Athens, Greece, 2007.

\bibitem
{elia.mitter.tac01} N.\ Elia and S.K.\ Mitter.
\newblock Stabilization of linear systems
with limited information.
\newblock {\it IEEE Trans.\ Automat.\ Control}, {\bf 46},  1384-1400, 2001.

\bibitem
{marino.tomei.I.tac93}
R.\ Marino and P.\ Tomei.
\newblock Global adaptive output-feedback
control of nonlinear systems, Part I.
\newblock {\em IEEE Trans.\ Auto.\ Contr.},
38, 17--32, 1993.

\bibitem
{marino.tomei.II.tac93}
R.\ Marino and P.\ Tomei.
\newblock Global adaptive output-feedback
control of nonlinear systems, Part II.
\newblock {\em IEEE Trans.\ Auto.\ Contr.},
38, 33--48, 1993.

\bibitem
{nair.et.al.proc.ieee07}
G.\ Nair, F.\ Fagnani, S.\ Zampieri, and R.J.\ Evans.
\newblock Feedback control under data rate constraints: An overview.
\newblock In {\em Proc.\ 2007 of the IEEE}, 95: 108--137, 2007.


\bibitem
{teel.praly.sicon95} A.R.\ Teel and L.\ Praly.
\newblock Tools for
semi-global stabilization by partial state and output feedback.
\newblock \emph{SIAM Journal on Control and Optimization}, 33, 1443--1488, 1995.




\end{thebibliography}
\end{document}